\documentstyle[11pt,epsfig]{article}
\begin{document}
%%%%%%%%%%%%%%%%%%%
\bigskip

\begin{center}
{\Large \bf Cutting Mutually Congruent Pieces from Convex Regions}\\
{\bf R Nandakumar}
{\normalsize (nandacumar@gmail.com)}\\
\vskip 0.4cm

{\bf ABSTRACT}
\end{center}
\textit{Question:} What is the shape of the 2D convex region $P$ from which, when 2 mutually congruent convex pieces with maximum possible area are cut out, the highest fraction of the area of $P$ is left over?\\

\textit{Partial Answer:} When $P$ is restricted to the set of all possible triangular shapes, our computational search gives an approximate upper bound of $5.6\%$ on the area fraction wasted when any triangle is given \textit{its} best (most area utilizing) partition into 2 convex pieces. We then give evidence for the general convex region which wastes the most area for its best convex 2-partition \textit{not} being a triangle and discuss some further generalizations of the above question.\\

\section{Introducing the Problem}

Two planar regions are \textit{congruent} if one can be made to perfectly coincide with the other by translation, rotation or reflection (flipping over).\\

\textbf{The Problem:} Which is the convex shape $P$ for which the largest fraction of its area gets left over ('wasted') under a partition of it into 2 congruent pieces with the largest possible area? \\

We call the fraction of the area of a given $P$ covered by 2 mutually congruent pieces cut from it as the 2-coverage of that partition. Any $P$ will have some such partition that maximizes this 2-coverage. We try to find $P$ which has the \textit{least maximum 2-coverage}. We consider only convex pieces; our approach is strongly experimental. To our knowledge, this problem was first stated in \cite{ref1} and \cite{ref2}.\\

\section{A Special Case - Triangles}

We first try to find the triangular shape with least maximum 2-coverage.\\

\textbf{Lemma 1:} There exist triangles such that their maximum 2-coverage with convex pieces necessarily less than a perfect 1.\\

\textbf{Proof:} Consider a scalene triangle $ABC$ with all angles acute and all sides of different length with $AB$ being the longest. If it can be perfect congruent partitioned into 2 convex pieces, it is easy to see that (1) side $AB$ has to be shared between the two pieces (2) the two pieces are separated by a straight line, say $l$ (3) line $l$ passes through a point on $AB$, say $P$ and also the opposite vertex $C$ - for if $l$ does not pass thru $C$, the triangle will be cut into a triangle and a quadrilateral, obviously non-congruent to one another. Now, $l$ intersects $AB$ at right angles - otherwise, one of the angles made at $P$ by $l$ will be obtuse and the piece without this angle will all angles acute and so the 2 pieces cannot be congruent. So, the only candidate 2-partition of $ABC$ is into two right triangles and they are not congruent since their hypotenuses are unequal - since all sides of $ABC$ are unequal. Further, we see that there is no line that splits $ABC$ into two convex congruent pieces (zero left over) such that the two pieces approach congruence to one another arbitrarily closely. $\diamond$\\

\textbf{Corollary:} There exist triangles for which the maximum 2-coverage with convex congruent pieces \textit{cannot} even be made arbitrarily close to 1.\\

Indeed, If two convex and mutually congruent tiles lie in the interior of triangle $ABC$ as described above such that they cover the interior arbitrarily closely, then, there is a straight line say $l$ that separates the tiles (because the tiles are convex). This $l$ divides $ABC$ itself into two convex pieces $P_1$ and $P_2$ and each of $P_1$ and $P_2$ contains a tile. Since the two tiles are congruent,$P_1$ and $P_2$ should be arbitrarily close to being congruent to each other, which is not possible both should arbitrarily approach triangles which cover $ABC$ perfectly and that is not possible by lemma 1.\\

\textbf{Inference:} From above, we surmise that there is a lower bound (a positive value less than 1.0) to the best convex congruent 2-coverage of triangles. We now describe an experimental approach to estimate this bound and the triangular shape that yields it. \\

Among the infinitely many ways in which any given triangle $T$ can be partitioned into 2 convex congruent pieces, we (mainly) consider only \textit{3 separate sets of partitions}:

\begin{enumerate}
\item \textit{Method 1:} using an angular bisector to cut 2 congruent triangular pieces from $T$ - a total of 3 such candidate partitions. In each partition, a small triangular bit goes waste.
\item \textit{Method 2:} with a suitable line parallel to an edge of $T$ as reference, cut 2 congruent quadrilaterals - this method yields a further 3 candidates. 2 triangular pieces are left over in each candidate partition.
\item \textit{Method 3:} partition into 2 congruent pentagons, 3 more candidates. Each candidate wastes 3 small triangles.
\end{enumerate}

Figure 1 shows one candidate of each type. For methods 1 and 2, the `cut lines' are shown dashed. Exception: line $A^\prime B^\prime$ in method 2 (partition into quads) is not a partitioning line but a reference line for the partition. The two quadrilateral pieces in method 2 are ${A Q^\prime Q A^\prime}$ and ${P B^\prime Q Q^\prime}$. In method 3, the 2 congruent pentagons cut from triangle $ABC$ are shown in orange and green colors.\\
%%%%%%%%%%%%%%%%%%%%%%%%%%%%%%%%%%%%%%%%%%%%%%%%%%%%%%%%%%%%%%%%%
\begin{figure}[t]
%begin{center}
\hskip 0cm
\centerline{\includegraphics[width=15cm, height=14cm, angle=0]{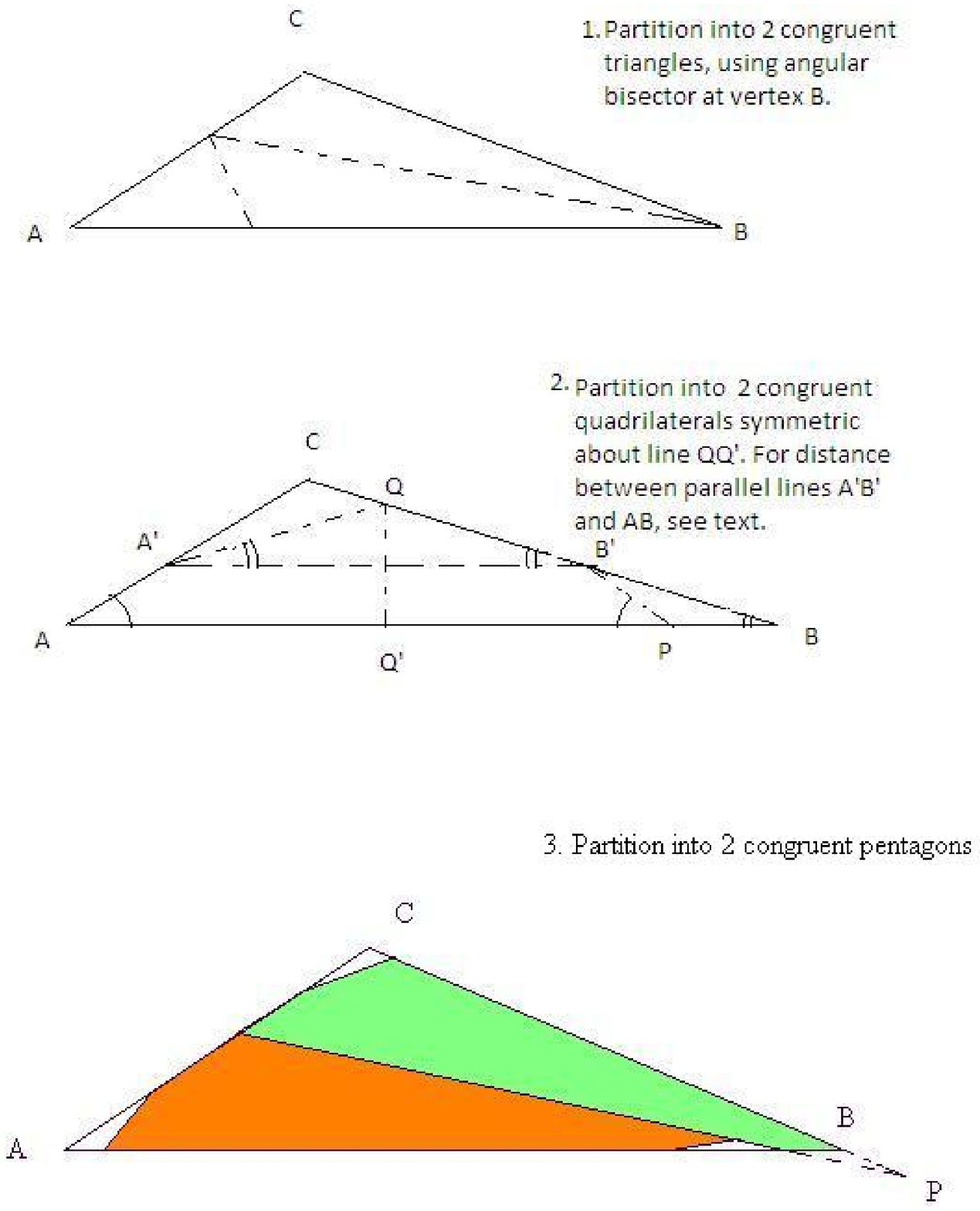}}
\caption{\label{} }
%\end{center}
\end{figure}
%%%%%%%%%%%%%%%%%%%%%%%%%%%%%%%%%%%%%%%%%%%%%%%%%%%%%%%%%%%%%%%%%

\textit{Description of Method 2:} if $\alpha$ and $\beta$ are the angles at base vertices $A$ and $B$ of the full triangle and $h$, the perpendicular distance of vertex $C$ from base $AB$, line $A^\prime B^\prime$ is parallel to the base $AB$ at a perpendicular distance $d$ above $AB$, given by:
\begin{center}
$d = h / (c+1)$ where $c = 2 \sin\alpha \cos\beta / \sin(\alpha+\beta)$.\\
\end{center}
The above expression gives the separation between the base and the reference line for which the coverage of the triangle by the 2 congruent quadrilaterals is the most. The derivation: Consider line $A^\prime B^\prime$ passing through $ABC$ at any given distance from $AB$ and parallel to it. If $A^\prime B^\prime$ is displaced by a small distance perpendicular to itself from this reference position, the area of one of the 2 triangular bits left out in the partition generated by $A^\prime B^\prime$ increases and the area of the other left out triangle decreases. The above value of $h$ is such that these two changes cancel - from this $h$, if $A^\prime B^\prime$ is parallel-displaced, the first order change in the 2-coverage of full triangle $ABC$ vanishes.\\

\textit{Short Description of Method 3:} The 2 congruent pentagons are calculated by an exhaustive search within the triangle. The 2 pentagons are obviously not mirror images of each other. We first find 2 congruent triangles, one of them lies entirely inside the  triangle $ABC$ and the other has a small portion projecting out. The point $P$ in figure 3 is this external vertex of one of these triangles. These 2 congruent triangles are then trimmed (where needed) and expanded (where possible) resulting in 2 congruent pentagons. There are 2 other candidate partitions of this type - in those cases, the outer point $P$ lies close to the other two vertices of the full triangle.\\

\pagebreak

\textbf{The Setup:}\\

All triangular \textit{shapes} can be generated by fixing the longest side (equivalently 2 of the vertices) of the triangle and only varying the position of third vertex within a finite region. We fix the longest side of the triangle to run from $(0, 0)$ (vertex $A$ in the earlier discussion) to $(10, 0)$ - the vertex $B$. The third (and only variable) vertex $C$ = $(x_3, y_3)$ varies inside a portion of a circle with radius 10 and centered at $B$ (10, 0) as shown in figure 2. It is easy to see that placing the third vertex within the yellow region exhausts all possible \textit{triangle shapes} up to reflections.\\
%%%%%%%%%%%%%%%%%%%%%%%%%%%%%%%%%%%%%%%%%%%%%%%%%%%%%%%%%%%%%%%%%
\begin{figure}[h]
\begin{center}
\hskip 0cm
\includegraphics[width=14.0cm, height=14cm, angle=0]{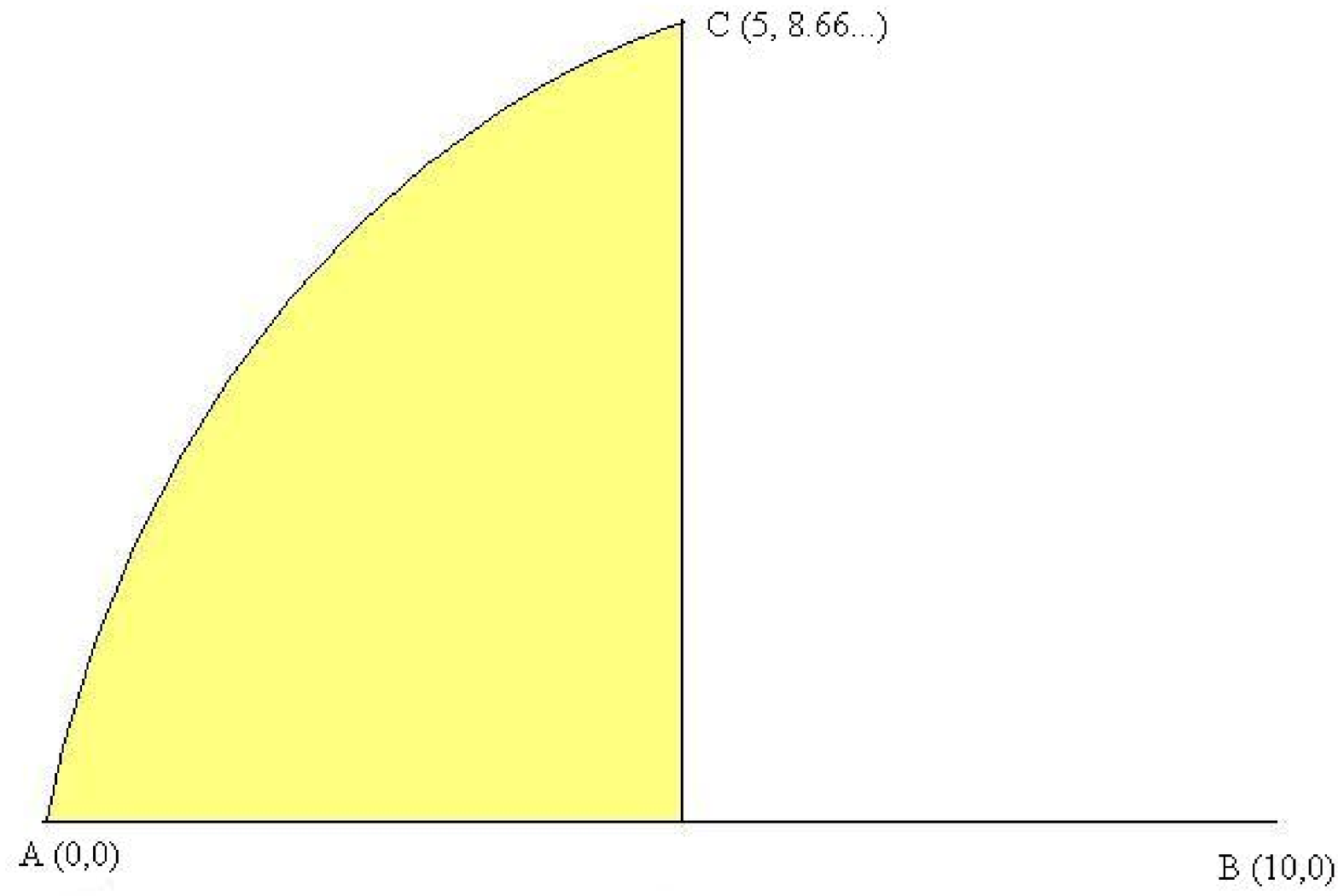}
\vskip -3cm
\caption{\label{} }
\end{center}
\end{figure}
%%%%%%%%%%%%%%%%%%%%%%%%%%%%%%%%%%%%%%%%%%%%%%%%%%%%%%%%%%%%%%%%%

For any position of $C$ $(x_3, y_3)$ in the yellow region in figure 2 (ie. for every shape of the triangle), we partition the resulting triangle $ABC$ using all candidates from all methods 1, 2, 3 and select the partition which gives the maximum 2-coverage of that triangle. Then, from maximum 2-coverages for all possible triangles, we select the minimum; the corresponding triangle is output.\\

\pagebreak

\textbf{Findings:}\\

For each candidate partition, we find that the 2-coverage of triangle $ABC$ for that partition as a function of the position of $C$  has a regular behavior with no multiple local maxima and minima as $C$ varies within its domain. So there is no threat of making qualitative errors if we restrict $C$ to a closely spaced grid of points.\\

1. For every triangle, if we consider only the 3 candidates from method 1 (using angular bisectors), the least maximum 2-coverage (equivalently, the maximum of the least wastage) is given a sliver (degenerate) triangle - which is also the triangle with the highest possible \textit{scalenity}; its maximum 2-coverage is $1 - 1/ \phi^3 =0.763.. $ ($\phi$ is the golden ratio). In our setup, the corresponding position of vertex $C$ = $(x_3,y_3)$ is $(3.82.., \delta)$ with $\delta$ tending to 0. However, this sliver has a partition by method 2 into congruent quadrilaterals resulting in max 2-coverage of almost 0.9 and indeed, it is not the triangle we are looking for.\\

2. If we only consider candidate partitions selected from methods 1 and 2, we find another sliver with sides tending to the ratio: $1 : 1/\sqrt2 : (1- 1/ \sqrt2)$ to have the highest least wastage. The best 2-partition of this shape has gives a 2-coverage of 0.8284...(a considerable increase from 0.763..) The corresponding Vertex $C$ $(x_3,y_3)$ is at $(2.92..., \delta)$ where $\delta$ tends to 0.\\

3. Finally, when we try all candidates from all the 3 methods for every triangular shape, the triangle with least maximum 2-coverage turns out to be `fat'. For $C$ at $(4.2, 6.7)$, the best partition of the resulting fat triangle gives a maximum 2-coverage of approximately 0.942 - only just under 6\% of the area of this triangle goes waste - and for every other triangle (got by varying $C$ in a lattice), the best 2-coverage is even higher than 0.942.. (and wastage, correspondingly less).\\

For this most wasteful triangle (which does not waste too much!), the best 2-coverage is given by 2 different candidate partitions - the pentagonal (method 3) partition with $B$ being the closest vertex to the external point $P$ and the partition into 2 congruent triangles given by  the bisectors of the angle at $C$. Both best congruent partitions are shown in Figure 3 - with the left out bits colored (not to scale).\\
%%%%%%%%%%%%%%%%%%%%%%%%%%%%%%%%%%%%%%%%%%%%%%%%%%%%%%%%%%%%%%%%%
\begin{figure}[t]
\begin{center}
\vskip -6cm
\includegraphics[width=10.0cm, height=18cm, angle=0]{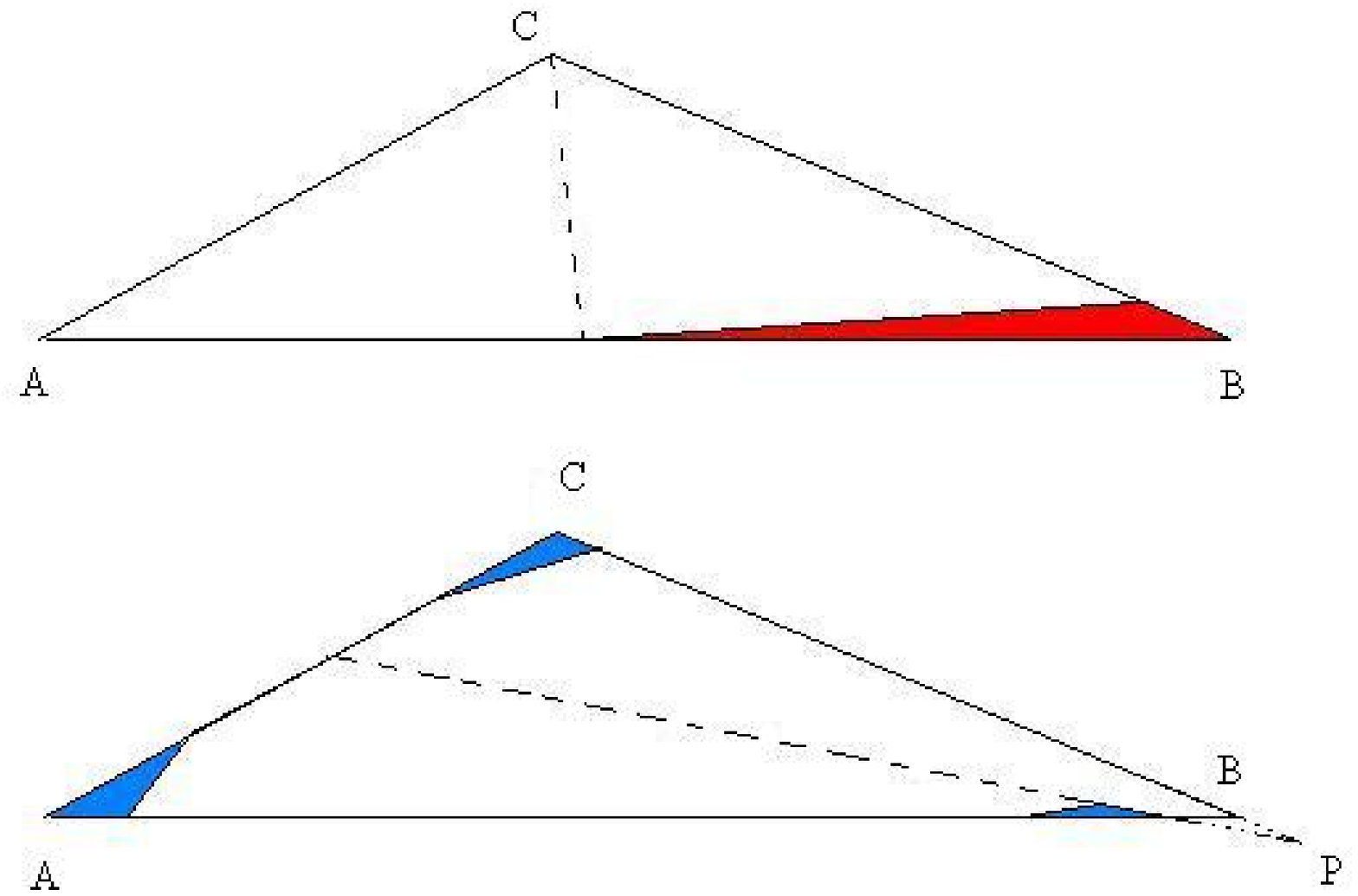}
\vskip -4cm
\caption{\label{} }
\end{center}
\end{figure}
%%%%%%%%%%%%%%%%%%%%%%%%%%%%%%%%%%%%%%%%%%%%%%%%%%%%%%%%%%%%%%%%%

The vertices $A$ and $C$ of the full triangle are not part of either congruent piece in the best pentagonal partition of this triangle; vertex $B$ = $(10, 0)$ is part of one of the pentagonal pieces.  On the other hand, vertices $A$ and $C$ are part of the best partition into 2 congruent triangles but vertex $B$ is left out by that partition.\\

\textit{Note 1}: We have a division of  the yellow zone above into 3 different regions: (1) positions of vertex $C$ = $(x_3, y_3)$ such that the congruent partition of $ABC$ with maximum 2-coverage (and least wastage) is given by the method 1 above (the one using angular bisectors, yielding 2 triangle pieces); (2) positions of $C$ where, the maximum 2-coverage is given by method 2 (quadrilateral pieces); (3) positions of $C$ where method 3 (pentagon pieces) gives the best 2-coverage.  These regions are separated by curves which converge on a 'triple point' which is approximately $(4.5, 5.3)$.\\

\textit{Note 2:} If we were to try still more partition schemes, for any given triangular shape, the maximum 2-coverage obviously cannot decrease. So the minimum among the max 2-coverages of all triangles can only increase from 0.942... As noted earlier, with only convex pieces, the least max 2-coverage cannot arbitrarily approach 1 for every triangle. So 0.942 ... is very close to the final answer. \\

\textit{Note 3:} `Near' the partition of a triangle into 2 pentagons (method 3), there are partitions into convex polygons with more sides (eg. it can be seen that the pentagonal pieces could be deformed into hexagons leaving out 4 tiny bits from the full triangle). Searching for these partitions and finding the best among them for each triangle could be computationally very expensive and unlikely to improve the bound substantially. However, for the most wasteful triangle, all such close variations on the best pentagonal partition necessarily leave out a small neighborhood of vertex $A$ - the point $(0, 0)$.\\

\section{Generalization}

We now consider the wider question: the general convex 2D shape that minimizes the maximum 2-coverage on congruent partitioned into 2 pieces of maximum area.\\

\textbf{Claim:} In the immediate neighborhood of the most wasteful triangle found by our search, we can find a convex shape with lower maximum 2-coverage than the triangle itself.\\

\textbf{Supporting Arguments:} Consider again, the most wasteful triangle with $C$ at (4.2, 6.7) and its 2 best partitions - one into congruent triangles using the angular bisector at vertex $C$ and another into congruent pentagons with the external corner point $P$ lying outside vertex $B$ (figure 3).\\

Let the most wasteful triangle (call this $T$) have area 1 in suitable units. Let a total area of $\alpha$ be covered by the 2 pieces in both the best partitions as in figure 3($\alpha$ is nearly .942). Trim $T$ slightly at both vertices $A$ and $B$ causing a loss of area of say, $\epsilon$ near each $A$ and $B$ resulting in another convex polygon, say $P$, with area $1-2\epsilon$. Since $P$ is in the immediate neighborhood of the most wasteful triangle, the partitions of it which maximize 2-coverage will be suitable slight deformations of the pieces of the best partitions of the triangle. It is easy to see that the total area of $P$ covered by the two congruent pieces under suitable deformations of either best partition of $T$ is $\alpha - 2 \epsilon$. So the 2-coverage for $P$ under \textit{its} best partitions is $(\alpha - 2\epsilon) / (1 - 2\epsilon)$. This is slightly less than $\alpha$, the maximum 2-coverage of $T$, the most wasteful triangle. Thus $P$ is a convex shape with lower 2-coverage than $T$.\\

We said above that some slight modifications of the pentagonal partition may slightly increase the 2-coverage for triangles resulting in a new most wasteful triangle, say $T'$, quite close to $T$. But these new partitions also leave out vertex A of the full triangle $T'$ and include B; so, the above trimming will give a convex polygon with lower maximum 2-coverage than $T'$.$\diamond$\\

\textbf{Conjecture:} If a convex polygon $Q$ is found with 3 optimal convex 2-partitions with coverage say $\beta$, such that it needs to be trimmed at 3 vertices for both pieces in \textit{all} the three best 2-partitions to lose $\epsilon$ area each, then, if we trim $Q$ at 3 vertices, its coverage changes from $\beta$ goes to $(\beta - 2\epsilon) / (1 - 3\epsilon)$. This new value of coverage is greater than  $\beta$ unlike the case of trimming the most wasteful triangle (in above proof) and such a $Q$ is at least a local minimum for the highest 2-coverage. This implies there could be a convex shape $R$ with $n$ different optimal 2-partitions (all of coverage $\gamma$, say) and which needs to be trimmed at $n$ vertices for all of these to be prevented. This means for $R$, the coverage \textit{increases} from $\gamma$ to $(\gamma-2\epsilon) / (1 - n\epsilon)$. Such a shape $R$, with $n$ tending to infinity (and hence $C^1$ smooth) could be the global minimum for maximum 2-coverage. This guess appears difficult to verify with our present knowledge.\\

\section{Conclusions}
The main question of the convex 2D shape that has the least maximum 2-coverage remains open although we are almost sure it is not a triangle and have guessed it would be a region with $C^1$ smooth boundary. We do not have much idea about number of pieces being larger than 2, even for partitioning triangles (except for some special cases) or about higher dimensions. Allowing the pieces to be non-convex could further increase the least maximum coverage.\\

\textbf{Claim:} We restate here a claim from \cite{ref2}: for any given number $N$, if any given convex polygonal region $P$ allows a zero-waste partition into $N$ non-convex congruent pieces each with finitely many sides, then $P$ also allows a partition into $N$ convex congruent pieces with zero wastage.\\

\textbf{Remarks}: Consider $N$=2. There exist convex shapes for $P$ which do not allow a perect (waste free) division into 2 pieces even with non-convex pieces. Indeed, such a partition has to divide the outer boundary of $P$ into exactly 2 connected curves and it is easy to see that with the same triangle example of lemma 1, these two connected curves from the boundary cannot not be mutually congruent. \textit{Guess:} however, it may be possible, given any convex polygon (not only triangles) and any $N$, to \textit{approach a perfect congruent partition arbitrarily closely} with non-convex pieces if the pieces could have arbitrarily many sides. This issue was briefly mentioned in \cite{ref2} where it was also speculated that in higher dimensions, the behavior could be qualitatively different.\\

\end{document}